\documentclass{article}
\usepackage{graphicx} 
\usepackage{listings}
\usepackage{hyperref}
\usepackage[mathletters]{ucs}
\usepackage[utf8x]{inputenc}
\usepackage{subfigure}
\usepackage{amsmath}
\usepackage{amssymb}
\usepackage{xfrac}
\usepackage{subfigure}
\title{Numerical Methods for Chaotic ODE}
\author{Robert M.~Corless\\
Dept. Computer Science\\
Western University}
\date{January 2025}
\newcommand{\e}{{\varepsilon}}

\begin{document}

\maketitle
\begin{abstract}
This paper explores backward error analysis for numerical solutions of ordinary differential equations, particularly focusing on chaotic systems.  Three approaches are examined: residual assessment, the method of modified equations, and shadowing. We investigate how these methods explain the success of numerical simulations in capturing the behavior of chaotic systems, even when facing issues like spurious chaos introduced by numerical methods or suppression of chaos by numerical methods.  Finally, we point out an open problem, namely to explain why the statistics of long orbits are usually correct, even though we do not have a theoretical guarantee why this should be so.
\end{abstract}

\section{Introduction}
Ordinary numerical methods work well for chaotic differential equations, even when they are not guaranteed by standard theory to do so, if some simple precautions are taken.  The simplest precaution to take is to use a variable-stepsize method and solve the problem more than once, with different error tolerances.  One will usually see---even if the computed trajectories are different from each other---that some statistics of the solution, for instance the dimension of the attracting set or the measures on the attracting set, are robust.  It's actually as simple as that.

To explain why this is so, however, is not so simple.  Indeed, parts of the explanation are not yet known.

Before we begin, I give some definitions.  The definition I give here of a ``chaotic system'' is nonstandard, but I claim it is equivalent to the standard one.

\textbf{Definition 1}: Let $f: \mathbb{R}^n \to \mathbb{R}^n$ be smooth, let $\e_1 >0$, $\e_2>0$, and $v_k:\mathbb{R} \to \mathbb{R}^n$ for $k=1$, $2$, be bounded for $t>0$; say $\|v_k(t)\|_\infty \le 1$.  A dynamical system
\begin{equation}\label{eq:firstorder}
  \dot{y} = f(y)
\end{equation}
is \textsl{chaotic} if its solutions are bounded and if typical solutions of the \textsl{persistently disturbed} dynamical systems
\begin{align}\label{eq:persistent}
  \dot{y}_1 &= f(y_1) + \e_1 v_1(t) \\
  \dot{y}_2 &= f(y_2) + \e_2 v_2(t)
\end{align}
with $y_1(0)=y_2(0)$ separate by $O(1)$ after time $T$ proportional to $\ln(1/\e)$ where $\e$ is the maximum of $\e_1$ and $\e_2$.  That is, $y_1(0)=y_2(0)$ but already by $t=T$, a short time later by virtue of the logarithmic dependence on $\e$, $y_1(t) - y_2(t) = O(1)$.

This definition sidesteps the usual difficulty with understanding how errors can grow exponentially while the solutions stay bounded; here we say that the solutions depart from one another after only a short time.  The reason this agrees with the usual definition is that if the system has a positive largest Lyapunov exponent, say $\lambda > 0$, then the difference $y_1(t)-y_2(t)$ will grow like $\exp(\lambda t)$.  But by the standard Gr\"obner--Alexeev nonlinear variation of constants formula~\cite{Hairer(1993)}, there exists a function $G$ such that
\begin{equation}\label{eq:GAVC}
  y_1(t) - y_2(t) = \e \int_{0}^{t} G(t,\tau) (v_1(\tau)-v_2(\tau))\,d\tau\>.
\end{equation}
This is $O(\e)$ near the beginning because $v_1$ and $v_2$ are $O(1)$; but once $t=T=\log(1/\e)/\lambda$ this will no longer be $O(\e)$ but rather $O(1)$.

One normally sees systems defined as being chaotic if the solutions with two different initial conditions depart from each other in this manner; but this is equivalent because in essence $G = \partial{y}/\partial{y_0}$, the derivative of the solution with respect to the initial condition.  We do it this way because the notion of \textsl{persistent disturbance} maps very nicely onto what numerical solutions do.

\subsection{An example: the Lorenz system}
I have claimed that numerical solution of chaotic dynamical systems is simple, if simple precautions are taken.
To demonstrate that this is so, we give a short script in Julia, which we take almost straight from the tutorial at \href{https://docs.sciml.ai/DiffEqDocs/stable/getting_started/#Defining-Parameterized-Functions}{the Getting Started Page}.  In it, we solve the Lorenz equations~\cite{Sparrow:1982}.
We should first ascertain that the Lorenz system is chaotic for the parameter values we'll use, but this is well-known.  Indeed, we can see the largest Lyapunov exponent being computed at \href{https://www.chebfun.org/examples/ode-nonlin/LyapunovExponents.html}{the Lyapunov exponent example on the Chebfun website} and with $T=25$ and $\e=10^{-9}$ so $\ln(1/\e) \approx 20.7$ (more about that in a moment) the result is $\lambda \approx 0.905 > 0$.  Notice that that website uses exactly the definition I described above.  It's written down in a few places (for example in the Proceedings of the Organic Mathematics Workshop from 1995, but sadly that's hard to find nowadays) but mostly people just use it because it's obviously equivalent to the standard definition.

Coming back to the example, I used the default method of the code described in~\cite{rackauckas2017differentialequations}, but changed the tolerances.
\begin{verbatim}
using DifferentialEquations
using Plots
function parameterized_lorenz!(du, u, p, t)
    x, y, z = u
    σ, ρ, β = p
    du[1] = dx = σ * (y - x)
    du[2] = dy = x * (ρ - z) - y
    du[3] = dz = x * y - β * z
end
u0 = [1.0, 0.0, 0.0]
tspan = (0.0, 5.0e1)
p = [10.0, 28.0, 8 / 3]
prob = ODEProblem(parameterized_lorenz!, u0, tspan, p)
sol = solve(prob, Tsit5(), reltol = 1e-8, abstol = 1e-8)
plot(sol, idxs = (1, 2, 3))
savefig("LorenzTsit588.png")
\end{verbatim}

\begin{figure}
    \centering
    \includegraphics[width=0.8\linewidth]{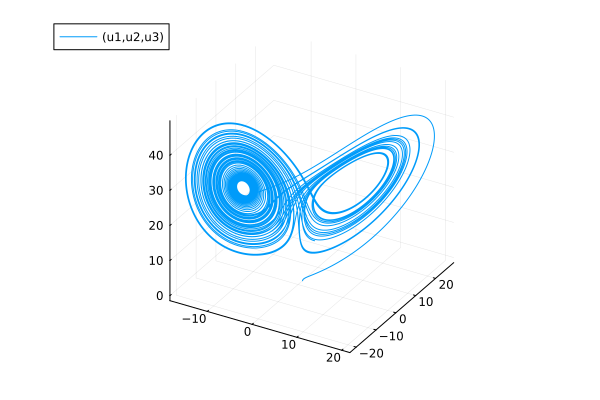}
    \caption{The computed solution to the Lorenz system on the interval $0 \le t \le 50$ with parameters $\sigma=10$, $\rho = 28$, and $\beta=8/3$, as computed by the Julia solver ``Tsit5'' with tolerances $10^{-8}$. With tolerances $10^{-9}$ the figure looks no different, to the casual eye.}
    \label{fig:LorenzTsit588}
\end{figure}

The plot appears in figure~\ref{fig:LorenzTsit588}.  It looks little different than any other plot of the solution to the Lorenz system computed by any other standard method.  There's a reason for that, and it is \textsl{not} because all the computed solutions agree with each other.  In fact they \textsl{don't} all agree with each other, in detail, because the Lorenz system is chaotic for these parameter values, as I said and as we all know.  In spite of that, the picture is familiar.  To confirm that, you can run Julia yourself, and change the method or the tolerances to anything you like, and you will see that the picture basically does not change.  But I will now demonstrate that the numerical values of the solution do indeed differ.
Issuing the command
\begin{lstlisting}
sol(50.0)
\end{lstlisting}
asks Julia to give you the values of $x$, $y$, and $z$ at $t=50$, and we get
\begin{lstlisting}
3-element Vector{Float64}:
 2.0436117819130932
 3.7394392065867637
 9.953318846465097 .
\end{lstlisting}
Doing the computation over again with the tolerances set instead to $10^{-9}$ gives instead when I asked for the new solution values
\begin{lstlisting}
3-element Vector{Float64}:
 -3.413851213342435
 -5.377552911950605
 17.822220963462385
\end{lstlisting}
and we see not a single decimal digit agrees in any of the components with what was computed before.  Indeed, you might get different results than the above if you run them on a different machine than the one I have run these on (a Microsoft Surface Pro) because even rounding errors can affect the results of a chaotic problem.
So why do the plots look the same?

\section{The answer is so simple it seems like cheating}
The reason is that both times the code gave us the exact solution, not to $\dot{Y} = F(Y)$ where $F$ stands for the three-dimensional Lorenz system and $Y=(x,y,z)$ is the vector of unknown functions, but rather $\dot{Y} = F(Y) + \e v(t)$, where $v(t)$ is a vector $O(1)$ in size and $\e$ is a ``small'' scalar.  That is, the numerical solution is the \textsl{exact} solution to a persistently disturbed Lorenz system!

That is the main point of this paper and this section.

Unfortunately, the $\e$ is not exactly $10^{-8}$ or $10^{-9}$ which is what we specified the tolerances to be, and it isn't even the case that the $\e$ when the tolerance is $10^{-9}$ is ten times smaller than it was when the tolerance was $10^{-8}$; if that were true then we would say the method had the property of ``tolerance proportionality,'' which is a nice property.  But what we do have is that the tighter the tolerance, the smaller the $\e$, generally speaking, and that turns out to be good enough.

Then we rely on a property of the Lorenz system that is well understood but not talked about much: its solution has quantities (called ``statistics'' in~\cite{Corless(1994)} and in~\cite{Corless(1994)c}) which are quite insensitive to persistent disturbances.

Think about it.  This must be true, or else the Lorenz equations would not be useful for any real purpose whatever except as a pure mathematical puzzle.  But in fact they provide quite a good model of, for instance, several analogue computers that were built to act as though the equations model them, and this in spite of the imperfections inherent in any physical object.  I typically summarize this idea as saying ``the model cannot be sensitive to trucks passing by the laboratory, or at least to fluctuations in the atmosphere of Jupiter which cause variations in its gravitational attraction.''  Real-world models \textsl{always} neglect small terms.  Some statistics of the model have to be insensitive to that neglect, or else the model is useless.

In fact, one of those robust statistics is the (fractal) dimension of the attracting set for the Lorenz system.  Lorenz himself proved that the attracting set was, well, attracting: he used a Lyapunov function to do so~\cite{lorenz1963deterministic}.  But more is true: the probability that the orbit will spend time in any portion of the attractor is also quite robust under perturbations of this sort.  To see that, you can do the computations yourself.

This is not only true of the Lorenz system, of course.  Any model of any real phenomenon must exhibit some kind of structural regularity with respect to perturbations (I use the phrase ``well-enough conditioned'' to indicate this because I am a numerical analyst by training, but a physicist might say that the model must be ``stable''.  Numerical analysts reserve that word for \textsl{algorithms} and say that an algorithm is numerically stable if it gives you the exact answer to a nearby question).

If the \textsl{trajectories} are insensitive to persistent disturbances, we say that the system is well-conditioned (and not chaotic).  Then all the other statistics will also be well-conditioned.  That's a sufficient condition.  But it is not a necessary condition.  Chaotic systems can have well-conditioned statistics, even though the trajectories are not.

Now we should look a little closer at the details, and at the $\e$.  How small is the $\e$ that the Julia system gave us, and how do we tell?  The method is to compute the \textsl{residual}, which is also called the ``defect'' in the older numerical analysis literature.  The residual is simply what's left over when we substitute a good interpolant to the ``skeleton'' of the solution back into the differential equation.

The word ``skeleton'' refers to the discrete set of points $(t_n, y_n)$ which are generated by the underlying numerical method.  The skeleton is what most people learn about when they first learn about numerical methods for ODE.  But modern codes provide interpolants automatically (for efficiency) so that really there is a continuously differentiable solution $y(t)$ available; if done well, this is mostly invisible to the user. The Julia codes are good codes, and they supply good interpolants so we don't have to worry about them.

\begin{figure}
    \centering
    \includegraphics[width=0.8\linewidth]{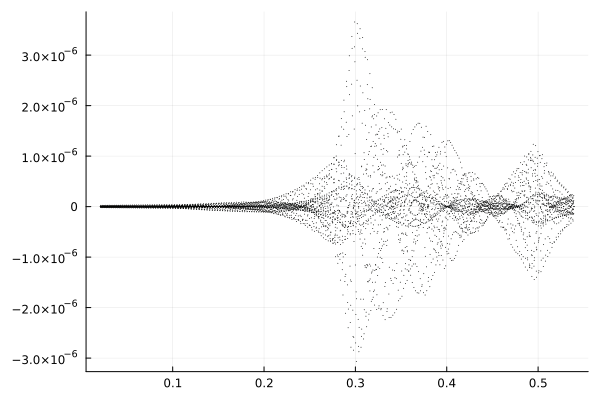}
    \caption{2025 samples of the residual in the computed solution to the Lorenz system on the interval $0 \le t \le 50$ with parameters $\sigma=10$, $\rho = 28$, and $\beta=8/3$, as computed by the Julia solver ``Tsit5'' with tolerances $10^{-10}$. }
    \label{fig:Tsit51010resid}
\end{figure}

The following script samples a few points ($N=800$) starting a few steps (15) away from the start up until a modest number of steps (240) after that.  It evaluates the solution at each of those $t$-values, and also the derivative of the interpolant to that solution at those values.  It then substitutes that derivative back into the original equation and computes what's left over: $r(t) = \dot{Y}-F(Y)$, which will be $r(t) = \e v(t)$ for the $O(1)$ vector $v(t)$.  So ``all we have to do'' is look at the size of $r(t)$ and that will tell us our $\e$.  First, we run this script with the output of the run with tolerances $10^{-8}$.
\begin{verbatim}
# Compute residual
N = 800
offset = 15
nsteps = 240
tsamp = Array(LinRange(sol.t[offset],sol.t[offset+nsteps],N));
res = Array{Float64}(undef,N,3);
σ, ρ, β = p;
for i=1:N
    (x,y,z) = sol(tsamp[i]);
    (dx,dy,dz) = sol( tsamp[i], Val{1});
    res[i,1] = dx - σ*(y-x);
    res[i,2] = dy - (x * (ρ - z) - y);
    res[i,3] = dz - (x * y - β * z);
end

# Plot the residual
plot( tsamp, res, seriestype=:scatter, markersize=0.1, legend=false )
\end{verbatim}
We see in the plot (not shown here) that the maximum value has magnitude about $1.5\times 10^{-4}$.  When we run the script again with tolerances ten times smaller, we find that the maximum value of the magnitude of the residual is about $2.2\times 10^{-5}$, which is smaller (as predicted) but not ten times smaller (as, sadly, was also predicted).  If we run it a third time, with tolerances $10^{-10}$, and this time with $N=2025$ because it makes a nice picture which we see in figure~\ref{fig:Tsit51010resid}, we find that the maximum residual is about $3.5\times 10^{-6}$.  Again, smaller, but not ten times smaller.

Of course, to be \textsl{sure} that the residual is always small, we would have to sample it from the start to the finish, and compute the maximum value attained.  This is not hard in this example; for a much larger system, it could be computationally expensive.  Each sample of the residual requires one ``function evaluation'' of $F(Y)$ and if the dimension of the system is large, then this will be costly.  But if you want to be \textsl{sure}, then you will pay the computational cost. To do that in this example, set the offset variable to be 1, set nsteps to be length(sol.t)-1, and then set N to be, say, 8, to take 8 samples of the residual on every step. Then use the \verb+maximum+ command to find the maximum residual.  You may also wish to compare it to the size of the variables $x$, $y$, and $z$, to compute a relative residual.

When we do this for tolerance $10^{-8}$, the maximum of $\sqrt{r_x^2 + r_y^2 + r_z^2}$ is $1.57\times 10^{-4}$.  For tolerance $10^{-9}$, the maximum is $2.36\times 10^{-5}$. For tolerance $10^{-10}$ the maximum is $3.67\times 10^{-6}$.  Computing this residual takes less time than plotting.

The idea that the numerical method has solved a different problem than the one you posed takes a little getting used to.  Wilkinson said that there is some evidence that people with mathematical training seem particularly unprepared to adopt the idea, and there is even some moral element to that judgement.  As I said, it almost seems like cheating.

\subsection{Does it always work?}
Almost always.  There would have to be something very special about the system, to have this idea fail: it would have to be very insensitive to the particular kind of persistent disturbance that numerical methods introduce, but sensitive to other kinds.  I have never seen an example of such.

It is \textsl{much} more common that you will fail if you solve a problem with a fixed step-size method and then when you check the size of the residual that you find it's gigantic, so what you have really done is solve a problem distant to the one you wanted.  In that case, you might get ``spurious chaos'' which is apparent chaos in a non-chaotic problem, or if you are using an implicit method for stiff problems then you might suppress existing chaos~\cite{corless1991numerical}.  See~\cite{Corless2016} for a more detailed introduction to what can go wrong.  See also~\cite{corless2019optimal} for what happens for stiff problems.

\subsection{Why doesn't $\e$ equal the tolerance?}
If you happen to be lucky enough to be using a defect-controlled code such as described in~\cite{Enright(2007)} then it \textsl{will} be (or, at least, almost always will be).  The downside of such a code is that it's a bit more expensive to run, and another downside is that it's not implemented in a nice free platform such as Julia. [ Apparently, the Julia codes \textsl{do} include one defect-controlled method, but I have never used it. ]

The real reason the residual has a different size than the tolerance is that controlling the size of the residual is asymptotically equivalent to controlling what is known as the ``local error per unit step'' which is more expensive in terms of computer time than controlling just the ``local error.'' The ``local error'' is a concept that makes sense only in the context of initial-value problems in ordinary differential equations, and if you are not going to build your own code it's really not helpful to explain it.  This ``local error'' is loosely related to the ``forward error'' or ``global error'' $y(t) - y_{\mathrm{ref}}(t)$ which, if the numerical method was going to converge to the reference solution as the tolerance went to zero, would be somewhat helpful.  Since that's what everyone tries to do, normally, then that's what gets implemented, even in the Julia package---which is so comprehensive and high-quality that I'm amazed every time I look at it.

For chaotic problems, though, convergence to the reference solution does not happen, except on algebraically short time intervals.  Trajectories of chaotic problems diverge exponentially (initially) from the reference trajectories, because their Lyapunov exponents are on average positive.

But we are lucky, in a sense.  Controlling the local error \textsl{indirectly} controls the residual~\cite{corless2013graduate}.  This is why I could say that tightening the tolerance would (generally) make the residual become smaller.

\subsection{Can we always compute the residual?}
Yes.  There is one difficulty, in that there is a different residual for every choice of interpolant for the problem, and some interpolants are better than others.  There is a concept of \textsl{optimal} interpolant, with the \textsl{minimal} backward error, and sometimes that is useful~\cite{Corless2018}, but usually the interpolant supplied with the code is good enough to estimate the residual accurately enough to assess the quality of the solution.  This is true for all the variable stepsize codes in~\cite{rackauckas2017differentialequations}, for example.

In the end, we will get the \textsl{exact} solution to
\begin{equation}\label{eq:AlteredEquation}
  \dot{Y}(t) = F( Y(t) ) + \e v(t)\>,
\end{equation}
where we will have computed $r(t) = \e v(t)$ and can plot it and see if it compares well with effects that the modeller neglected in setting up the equations $\dot{Y} = F(Y)$ in the first place.  We already know that the trajectories are exponentially sensitive to such changes, but what we will want to estimate is how insensitive the physically relevant statistics (eg dimension, measure on the attractor) are.  That's up to you, the modeller.

For the vast majority of your numerical solutions of chaotic problems, this process will be satisfactory.  But, sometimes, there are complications.
\section{Correlations with the solution}
We all know that if we solve $\ddot{y} + y = 0$ we get harmonic solutions; by shifting the origin if necessary we can say that $y(t) = A \cos(t -\phi)$ for some constants $A$ and $\phi$.  If we now persistently disturb this problem by a small amount $\e v(t)$, so we solve $\ddot{y} + y = \e v(t)$ instead, then for the vast majority of perturbation functions $v(t)$ the effect remains small (this problem is not chaotic).  But if we happen to perturb it at a \textsl{resonant} frequency, so $v(t)$ contains a term $\sin(t)$ or $\cos(t)$, then the solutions diverge linearly as $t \to \infty$.  This is called ``secular growth'' and really it's not a big problem, at least when compared to exponential growth as occurs in chaotic problems.  But it highlights that if the perturbation is correlated with the solution, then the effects of the perturbation may not be trivial.

Indeed the whole science of ``chaotic control'' is based on this, but in reverse: by artfully choosing the correlation of $v(t)$ with $y(t)$ you can make the trajectory do whatever you want! See e.g.~\cite{liao2006chaos}.

Numerical methods can introduce correlations of their errors with the solution.  Indeed most error control methods explicitly feed the solutions back into the stepsize to take next, and this is usually a feature, not a problem, of the method.

But if you are concerned, you might try to fit $\e v(t)$ to terms containing $Y(t)$ or $\dot{Y}(t)$, either experimentally (this is not a terrible idea). Working theoretically, if you know that the code is using a method that is amenable to analysis, you might be able to deduce the correlations ahead of time.

Suppose for instance that you are using a fixed step-size explicit Euler method with time step $h$, because your computer has lots of processors, your problem is very large, very stiff\footnote{Yes, it's possible for a problem to be both chaotic and stiff.}, and your simulation is never allowed to return a ``did not converge'' message and stop unexpectedly\footnote{In this case, explicit Euler is about as efficient as it gets, even though it's not very efficient.}.  In that case, the ``method of modified equations'' says that your code is more nearly solving
\begin{equation}\label{eq:modeuler}
  \dot{Y} = F(Y) + h J_F(Y) F(Y) + O(h^2) = (I + hJ_F(Y))F(Y) + O(h^2)\>.
\end{equation}
Now you could compute the residual in \textsl{this} equation (using, say, cubic Hermite interpolation) and identify that $O(h^2)$ term as your $\e v(t)$. We would expect $\e$ to be $O(h^2)$ then and not $O(h)$. At least, you could sample the residual in a few strategic places. The Jacobian matrix $J_F(Y)$ will be something that you can examine as the solution progresses so you can see if the term containing it is going to cause difficulty of interpretation.

For an introduction to the method of modified equations, see~\cite{corless2013graduate}.  For a more thorough treatment, see~\cite{Griffiths(1986)}.

\section{Shadowing}
A more reasoned objection to the previous analysis than saying ``it seems like cheating'' is to point out that in many models, the physical laws governing them are well-studied and have been proved to be valid experimentally to many decimal places.  Both general relativity and QED for example are known to give correct predictions to about eleven significant figures, and astronomical data is often even more precise.  For instance, one frequently thinks of gravity as being a wonderfully precise law, whereas it's the data to go into the equations that we don't know very well.

In this case, it can \textsl{still} make sense to assume that the data is really exact, and want to know what the laws predict.  One place where people will (sometimes grudgingly) accept less-than-perfect knowledge of the data is in the \textsl{initial conditions}.  Perhaps we don't know the exact positions of all the planets at the start of our solar system simulation (but we know all their masses exactly, and they're all exact point masses---sure, let's go with that).  Can we take the computed trajectory and find some exact reference trajectory that starts at a \textsl{slightly} different set of initial conditions but matches our computed trajectory for all time?  In that case, we say that our computed trajectory is \textsl{shadowed} by an exact trajectory.

It is generically true that locally-accurate numerical solutions are  shadowed~\cite{Pilyugin(1999)} for long times.  What does it mean, ``generic?''  Well, there is a technical meaning, which is complicated, because there are so many dynamical systems that it doesn't make sense to talk about ``measure'' or ``probability'' for them.  People interpret that technical meaning to be ``most'' or something like ``almost all.''  That helps, right?

Well, if it makes you feel better.  But the issue is that we will never know if any given dynamical system is ``generic'' or not.  It most likely is, but we won't usually know.  And we do know of some that are not.

Can we tell computationally if our computed trajectory is shadowed by a reference trajectory?  Yes, sometimes.  See~\cite{hayes2007fast} for the state of the art, and~\cite{hairer1997life} for a slightly different view.  It's more expensive than just computing a residual.  It starts with that, and then uses iterative refinement (something like Newton's method) globally.  It can be worth it, if you know the laws governing your model very well, but don't know your initial conditions precisely.  Indeed that's not uncommon.

\section{Hamiltonian Problems and Other Constrained Problems}
It gets more interesting when the ODE system is Hamiltonian or has other constraints and one must use methods that must satisfy the constraints, which are sometimes implicit.  One could write books about such problems: for instance, \cite{SanzSerna(1994)} on Hamiltonian problems; \cite{Hairer(2006)} for more general geometric integration (see also~\cite{crouch1993numerical}); and \cite{kunkel2006differential} for DAE (differential algebraic equations).  There is also the interesting new paper~\cite{kunkel2023discretization} which studies DAE with symmetries.

For such problems one is \textsl{much} better off with a method that preserves the symmetries or constraints as well as possible, and insists that whatever residual there is lie in the manifold specified by the constraints.  This is a kind of \textsl{structured} backward error analysis.

A simple example may suffice.  See also the discussion in~\cite{corless2013graduate}.

Consider the Henon-Heiles model.
This set of model equations was originally motivated by the motion of stars in a galaxy; however, the gravitational potential was chosen for simplicity, not detailed realism.  It is now used as an interesting example in several texts, for example~\cite[p.~187]{Bender(1978)} and~\cite{Hairer(2006)}.
The Hamiltonian is
\begin{align}
H = \frac{1}{2}\left( p_1^2 + p_2^2 + q_1^2 + q_2^2 \right) + q_1^2q_2 - \frac{1}{3} q_2^3\>.
\end{align}
We use the initial conditions given in~\cite{channell1990symplectic} at first: $p_1(0) = p_2(0) = q_1(0) = q_2(0) = 0.12$.  The equations of motion are $\dot{p} = -H_q$ and $\dot{q} = H_p$, or
\begin{align}
\dot p_1 &= -q_1 -2q_1q_2 \nonumber\\
\dot p_2 &= -q_2 - q_1^2 + q_2^2 \nonumber\\
\dot q_1 &= p_1 \nonumber\\
\dot q_2 &= p_2 \label{eq:HenonHeiles}
\end{align}
and we take a long(ish) time span of integration: $0 \le t \le 10^5$.  The constant-coefficient Hamiltonian is supposed to be conserved:
\begin{equation}\label{eq:conserveH}
  \frac{dH}{dt} = H_p \dot{p} + H_q \dot{q} = -H_p \cdot H_q + H_q \cdot H_p = 0\>.
\end{equation}

It turns out that this is already an example in the Julia tutorial.  See
\href{https://docs.sciml.ai/DiffEqDocs/stable/examples/classical_physics/#Henon-Heiles-System}{Henon-Heiles-System}.
That example demonstrates fairly convincingly that a symplectic method does a better job of conserving the energy, with less effort, than a non-symplectic method.  That example does point out that integrating at tight tolerances does a respectable job, though; it just takes more computing effort.

In~\cite{corless2013graduate} the ``leapfrog'' or St{\o}rmer--Verlet method is used, with fixed step-size $h$, as follows:
Following~\cite{Preto(1999)class} and also~\cite{Hairer(2006)}, we mean the `drift-kick-drift' version of leapfrog:
\begin{align}
\vec{q}^{(1)} = \vec{q}_n + \frac{h}{2} H_p(\vec{p}_n,\vec{q}_n) \enskip&\phantom{=}\>, \phantom{=}\vec{p}^{(1)} = \vec{p}_n\nonumber\\
\vec{q}^{(2)} = \vec{q}^{(1)} \enskip &\phantom{=}\>,\phantom{=}
\vec{p}^{(2)} = \vec{p}^{(1)}-h H_q(\vec{p}^{(1)},\vec{q}^{(1)}) \nonumber\\
\vec{q}_{n+1} = \vec{q}^{(2)} + \frac{h}{2}H_p(\vec{p}^{(2)},\vec{q}^{(2)}) \enskip &\phantom{=}\>,\phantom{=}\vec{p}_{n+1} = \vec{p}^{(2)}\>. \label{eq:DKD}
\end{align}
That is, we first drift with constant momentum for half a step, then we kick the system with the potential, and then we drift another half a step.

For second-order problems where the Hamiltonian is separable (as here) and furthermore $T(p) = \sfrac{(p_1^2 + p_2^2 + \cdots p_N^2)}{2}$, the equations have $\dot q_i = p_i$ and so are equivalent to a system of second-order ODEs for the $q_i$. Then, this method is equivalent to
\begin{align}
\vec{p}_{n+1} &= \vec{p}_n - h U_q(\vec{Q}_n)  \nonumber\\
\vec{Q}_{n+1} &= \vec{Q}_n + h  \vec{p}_{n+1}\>,
\label{eq:leapfrog}
\end{align}
as you can see by putting $\vec{Q}_n = \vec{q}_{n} + \sfrac{h \vec{p}_n}{2}$.
In either formulation, the method only requires one evaluation of the force term $U_q$ per step, just as Euler's method does.  For this class of problems, however, the method is second order, and \emph{symplectic} in that it gives samples at equal stepsizes $t_k = t_0 + kh$ of vector functions $\tilde{\vec{p}}(t)$ and $\tilde{\vec{q}}(t)$ that (nearly) satisfy the equations of motion from a \emph{perturbed} Hamiltonian system
\begin{align}
\tilde{H} = H(\vec{p},\vec{q}) + h^2 H_2(\vec{p},\vec{q}) + h^4 H_4(\vec{p},\vec{q}) + \cdots\>.
\end{align}
Carrying out the computation explicitly, we have (omitting terms that are $O(h^6)$)
\begin{align}
\dot p_1 &= -q_1 -2q_1q_2 +\frac{1}{6} K_{2} h^{2}+\frac{1}{60} K_{1} h^{4} \nonumber\\
\dot p_2 &= -q_2 - q_1^2 + q_2^2 +\frac{1}{12} K_{5} h^{2}-\frac{1}{120} K_{4} h^{4}\nonumber\\
\dot q_1 &= p_1 -\frac{1}{12} K_{8} h^{2}-\frac{1}{120} K_{7} h^{4}\nonumber\\
\dot q_2 &= p_2 -\frac{1}{12} K_{11} h^{2}-\frac{1}{120} K_{10} h^{4} \label{eq:HenonHeilesMod}
\end{align}
where
\begin{align}\label{eq:modterms}
K_{1} =&\>
-2 q_1 +5 p_1  p_2 -24 q_1  q_2 -20 q_1  q_2 ^{2}-24 q_2  q_1 ^{3}-8 q_1  q_2 ^{3}+5 q_1  p_1 ^{2}\nonumber\\
&\>{}-q_1  p_2 ^{2}-20 q_1 ^{3}+6 q_2  p_1  p_2  \nonumber\\
K_{2} =&\>
p_1  p_2 -2 q_1 ^{3}-2 q_1  q_2 ^{2}-6 q_1  q_2 -q_1 \nonumber\\
K_{4} =&\>
2 q_2  p_1 ^{2}-10 q_2  p_2 ^{2}+24 q_1 ^{2} q_2 ^{2}+40 q_2 ^{3}+12 q_1 ^{4}-20 q_2 ^{4}+40 q_1 ^{2} q_2 -12 q_1  p_1  p_2 \nonumber\\
&\>{}+24 q_1 ^{2}-24 q_2 ^{2}-5 p_1 ^{2}+5 p_2 ^{2}+4 q_2 \nonumber\\
K_{5} =&\>
p_1 ^{2}-p_2 ^{2}-6 q_1 ^{2}+6 q_2 ^{2}-2 q_2 -4 q_1 ^{2} q_2 -4 q_2 ^{3}\nonumber\\
K_{7} =&\>
12 q_2  q_1  p_2 -2 p_1  q_2 ^{2}+10 p_1  q_2 +10 q_1 ^{2} p_1 +10 q_1  p_2 +p_1 \nonumber\\
K_{8} =&\>
2 q_1  p_2 +2 p_1  q_2 +p_1 \nonumber\\
K_{10} =&\>
-2 q_1 ^{2} p_2 +12 q_1  p_1  q_2 +10 q_1  p_1 +10 q_2 ^{2} p_2 -10 q_2  p_2 +p_2 \nonumber\\
K_{11} &=\>
2 q_1  p_1 -2 q_2  p_2 +p_2
\end{align}
Integrating those, we find
\begin{align}\label{eq:H2}
  H_2 =&\>
-\frac{1}{12} p_{1}^{2} q_{2} -\frac{1}{24} p_{1}^{2}-\frac{1}{6} q_{1} p_{2} p_{1} +\frac{1}{12} p_{2}^{2} q_{2} -\frac{1}{24} p_{2}^{2}+\frac{1}{12} q_{1}^{4}+\frac{1}{6} q_{1}^{2} q_{2}^{2}+\frac{1}{12} q_{2}^{4}\nonumber\\
&\>{}+\frac{1}{2} q_{1}^{2} q_{2} -\frac{1}{6} q_{2}^{3}+\frac{1}{12} q_{1}^{2}+\frac{1}{12} q_{2}^{2}
\end{align}
and
\begin{align}\label{eq:H4}
  H_4 =&\>
\frac{1}{10} q_{1}^{4} q_{2} +\frac{1}{15} q_{1}^{2} q_{2}^{3}-\frac{1}{30} q_{2}^{5}-\frac{1}{24} p_{1}^{2} q_{1}^{2}+\frac{1}{120} p_{1}^{2} q_{2}^{2}-\frac{1}{10} p_{1} q_{2} q_{1} p_{2} +\frac{1}{120} p_{2}^{2} q_{1}^{2}\nonumber\\
&\>{}-\frac{1}{24} p_{2}^{2} q_{2}^{2}+\frac{1}{12} q_{1}^{4}+\frac{1}{6} q_{1}^{2} q_{2}^{2}+\frac{1}{12} q_{2}^{4}-\frac{1}{24} p_{1}^{2} q_{2} -\frac{1}{12} q_{1} p_{2} p_{1} \nonumber\\
&\>{}+\frac{1}{24} p_{2}^{2} q_{2} +\frac{1}{5} q_{1}^{2} q_{2} -\frac{1}{15} q_{2}^{3}-\frac{1}{240} p_{1}^{2}-\frac{1}{240} p_{2}^{2}+\frac{1}{60} q_{1}^{2}+\frac{1}{60} q_{2}^{2}\>.
\end{align}
That is, we can (spectrally nearly) interpolate the fixed-stepsize skeleton by using the solution of a nearby Hamiltonian problem.  The smaller $h$ is, the more work we have to do, but the nearer the perturbed Hamiltonian is to the one we wanted.  See~\cite{Calvo(1994)} as well as the previously-mentioned references for how to compute these perturbed Hamiltonians in general, but note that the method of modified equations as discussed in~\cite{corless2013graduate} and again in~\cite{CorlessFillion2025SIAM} will work for particular examples and methods.

For a variety of reasons, this is a more satisfactory backward error analysis than the simple method of computing the residual by a polynomial interpolant.  First, the perturbation of the problem is autonomous (if the original problem is autonomous---symplectic methods work for some time-dependent Hamiltonian problems as well). This means that the perturbed differential equation has the same dimension as the original.  Second, many physical perturbations of (for example, computational astronomy problems) are themselves Hamiltonian---think of neglecting the influence of other planets, for example.  Some, of course, are not, such as tidal friction or minor collisions, but these may be smaller than the Hamiltonian perturbations.  Third, preservation of symplecity (most people say ``symplecticness'') may better preserve certain statistical measures.  Given the chaotic nature of many Hamiltonian systems, of course, there is no hope of ensuring small global forward error in trajectories and so accurate statistics are all that can be computed.

Either way, one can look for reference trajectories that shadow the computed trajectory.  It's easier to find them if the method is a symplectic method.

\begin{figure}
  \centering
  \subfigure[$h=1.175$]{\includegraphics[width=0.45\textwidth]{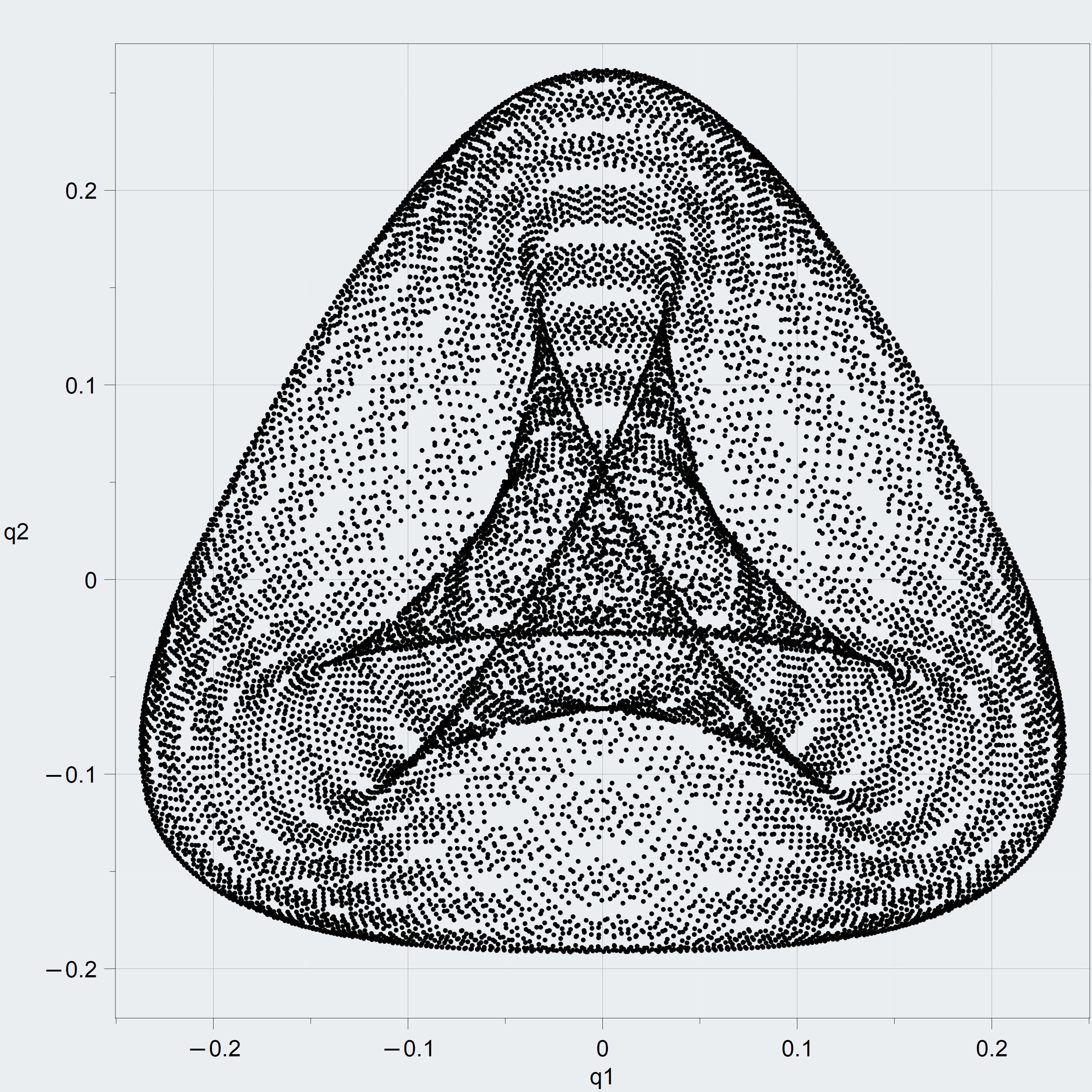}}
  \subfigure[$h=1.18$]{\includegraphics[width=0.45\textwidth]{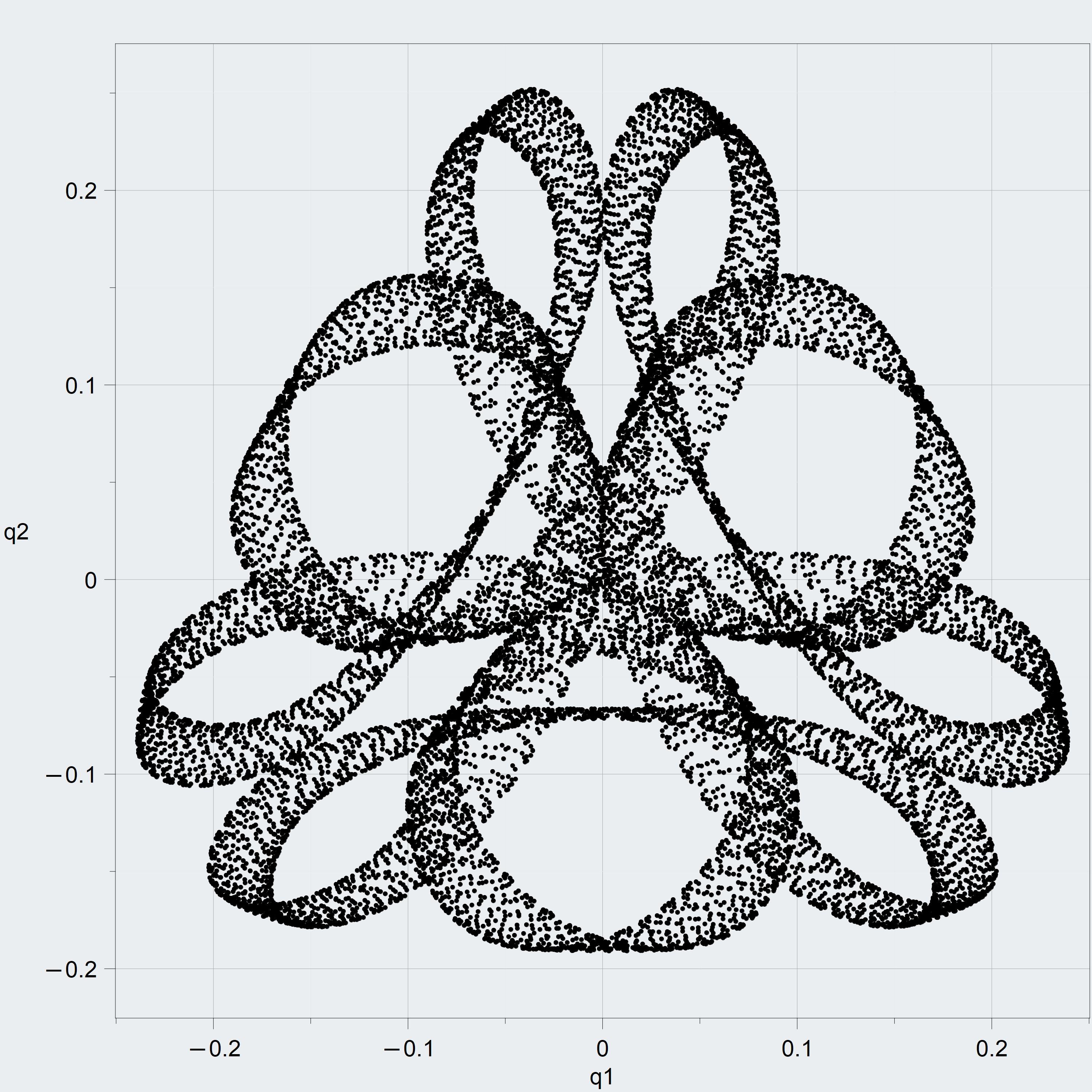}}\\
  \subfigure[$h=79/64$]{\includegraphics[width=0.45\textwidth]{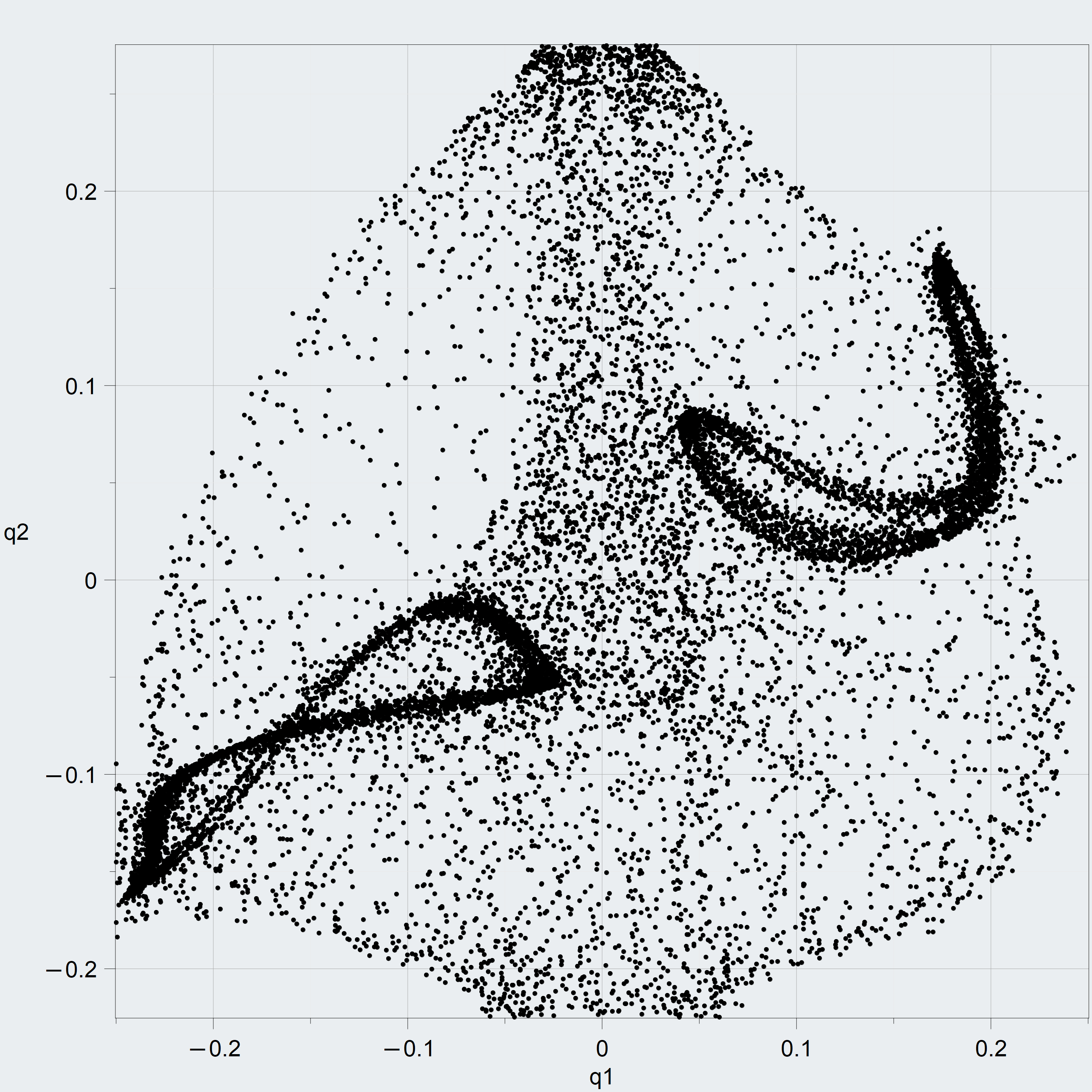}}
  \subfigure[$h=81/64$]{\includegraphics[width=0.45\textwidth]{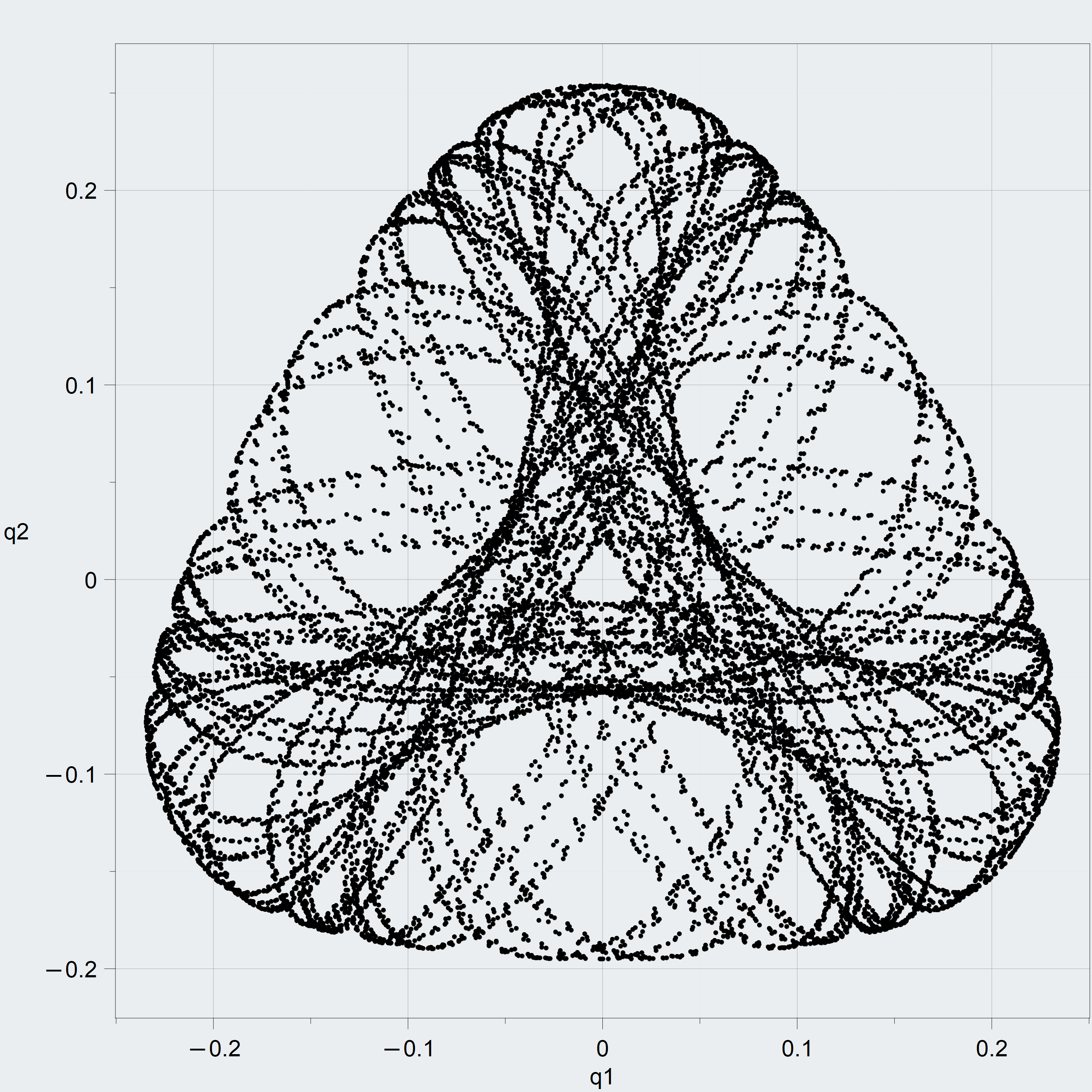}}
  \caption{Solving the Henon--Heiles equations with the given initial conditions using the St{\o}rmer--Verlet method with $h=1.175$ (top left) and $h=1.18$ (top right) for $N=16,000$ steps. Both have bounded state and nearly-constant energy, though resonance effects make the plots look quite different. Intermittent spurious chaos ensues with $h$ slightly larger, say by $h=79/64 = 1.234375$; by $h=81/64=1.265625$ stability has apparently returned.  In the graph, $q_2$ is plotted against $q_1$. For slightly larger $h$, all structure vanishes.  The energy of the solution remains reasonably constant over this time interval for three out of the four graphs. The energy starts varying noticeably when the trajectory is chaotic, so the spurious chaos is detectable.}\label{fig:HenonHeilesState}
\end{figure}

If we take $N=16000$ steps with $h=1.175$, which is near to where stability breaks down and the method blows up, it is still true that $H=H_0 + h^2 H_2 + h^4 H_4$ is more nearly constant than either $H_0$ or $H_0 + h^2H_2$ are, even though $h$ is larger than $1$.  At that $h$, the departure from the initial value is no more than $0.006$ even for $H_0$, and is smaller for $H_0 + h^2 H_2$, and smaller yet when the fourth order term is included. The state has the curious appearance shown in the top left-hand side of figure~\ref{fig:HenonHeilesState}, mostly owing to sampling effects and resonance.  Taking $h=1.18$, slightly larger, gives a completely different picture, for similar reasons; taking $h=79/64=1.234375$ (larger yet) we see apparent chaos.  This chaos is spurious, and the energy variations (not shown here) demonstrate this. The initial energy is about $0.034$ and the variations with $h=79/64$ are sometimes twice that.  Taking $h$ just a bit larger, $h=81/64 = 1.265625$, stability of a sort returns, and the energy variations are then no more than $0.001$ for $H_0 + h^2 H_2 + h^4 H_4$, larger for just $H_0$ and for $H_0 + h^2H_2$, which are about $0.009$ and $0.003$, respectively.  This intermittency is a well-known feature of many chaotic systems.

The main points of this example are first that even a good fixed-step numerical method can introduce spurious chaos, which is well-known, and second that measuring the energy makes the spurious chaos detectable.

\section{Is the shadowing orbit ``typical?''}
Here is where we run into some open problems.  The difficulty is most easily demonstrated with a \textsl{discrete} dynamical system, either the Gauss map, the logistic map, or the Bernoulli shift. The same thing happens with continuous dynamical systems, if they are simulated long enough.  I first ran into this in~\cite{Corless(1992)continued}, with the Gauss map, which is defined to be
\begin{equation}\label{eq:Gauss}
  G(x) = \mathrm{frac}( 1/x )
\end{equation}
if $x \ne 0$, and $G(x)=0$ if $x=0$.  This maps the half-open interval $[0,1)$ onto itself and the iteration $x_{n+1} = G(x_n)$ generates the continued fraction of $x_0$ as a byproduct.  I proved that the computed orbits, in IEEE arithmetic, were \textsl{uniformly} shadowed by a reference orbit.  All of the orbits are shadowed, with an orbit starting no more than $4$ times the unit roundoff away from $x_0$.  This is as strong a shadowing result as you can get.  Every computed orbit is shadowed by a true reference orbit, for all time.

But every single one of those computed orbits is \textsl{periodic} because the set of floating-point numbers is finite.  That means that the statistics of the orbits \textsl{must} be wrong, and \textsl{cannot} give the correct probabilities, because almost all true orbits are aperiodic and dense in the interval $[0,1)$.

The same is true for the logistic map $x_{n+1} = 4x_n(1-x_n)$, and for the Bernoulli shift map $x_{n+1} = 2x_n$ mod $1$.

\begin{figure}
    \centering
    \subfigure[fixed point]{\includegraphics[width=0.4375\textwidth]{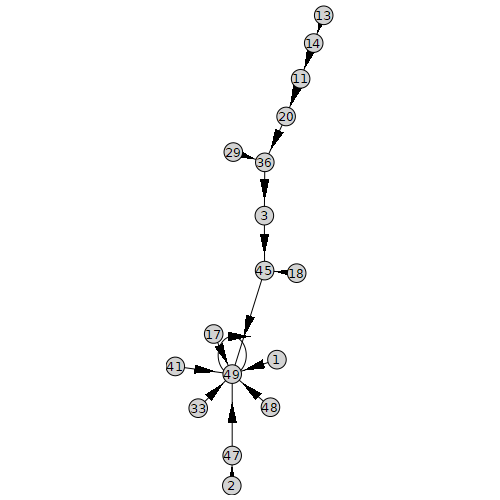}}
    \subfigure[two-cycle]{\includegraphics[width=0.4375\textwidth]{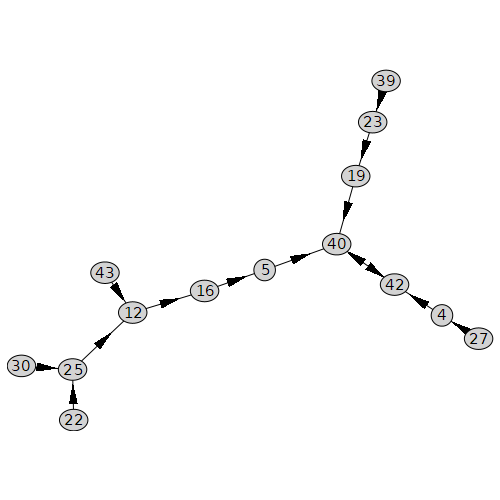}}\\
    \subfigure[another two-cycle]{\includegraphics[width=0.4375\textwidth]{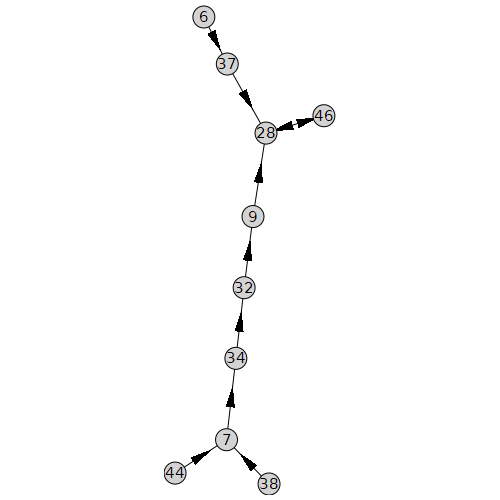}}
    \subfigure[three-cycle]{\includegraphics[width=0.4375\textwidth]{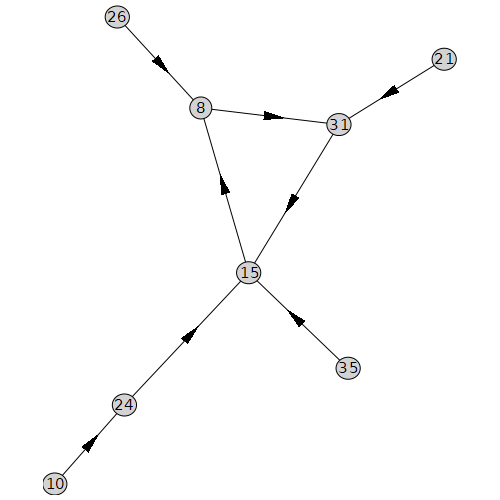}}
    \caption{All the possible orbits of $G(x)$ in quarter-precision IEEE floating-point.}
    \label{fig:ComparisonPlot}
\end{figure}

That may seem strange, but to demonstrate it in a simple way let's simulate the action of $G$ in $8$-bit IEEE Standard floating-point arithmetic (again I used Julia, but I used the GraphTheory package in Maple to draw the graphs).  I used $8$ bits, which is ridiculously low precision, to make the orbits explainable\footnote{But to my shock and horror, $8$-bit floats are actually used in machine learning.  Even worse, so are \textsl{four} bit floats, or three bit integers! Yet the results are claimed to be good. Possibly for good reasons, which we don't know yet.}.  There are only $49$ of these $8$-bit floats in the interval $[0,1]$, and the Gauss map takes each of these $49$ numbers to one of the others.  By direct simulation we can draw a graph on that set, and locate all the orbits.  We display them in figure~\ref{fig:ComparisonPlot}.

Half-precision has more points (15361) in $[0,1]$, more orbits, and longer cycles and longer transients.  Single-precision has yet more, and longer.  Double precision longer yet, and quad precision even better.  Yet all are finite, so this quarter-precision picture is not qualitatively wrong.  There are always a finite number of components, and at most a finite longest transient and a finite longest cycle. Moreover, the length of the longest cycle and the longest transient scales as the \textsl{square root} of the number of available floating-point numbers.  Even so, the statistics \textsl{aren't that bad}.

To repeat, the statistics of the orbits are not terrible approximations.  This is not fully explained, although some things are known~\cite{Boyarsky(1997)}: ``computers like Lebesgue measure'' in part because floating-point numbers are denser around $0$ than they are elsewhere.  See~\cite{Tupper2009} for an exploration of this in molecular dynamics simulations.  Again, short-term accuracy seems to give accurate statistics for long-term simulations, but nobody has a full explanation why.

Analysis of successful methods is not popular: ``if it ain't broke, don't fix it.''  The first method of this paper is a method that allows you to be \textsl{sure} that your results are good (you could even do it using interval arithmetic~\cite{Corless(1992)} if you wanted), and to detect the very rare cases when your results will not be good.  But, since most of the time it will be fine even if no check is made, most of the time people don't bother to check.  You will probably get away with it if you use a good numerical method and at least vary the tolerances occasionally.  But if money or lives are on the line, perhaps you had better check.  For further tools for analysis, consult~\cite{stuart1998dynamical}.

\section*{Acknowledgements}
This work supported by NSERC RGPIN-2020-06438.  The support of the Rotman Institute of Philosophy is also gratefully acknowledged.


\begin{thebibliography}{10}

\bibitem{Bender(1978)}
Carl~M. Bender and Steven~A. Orszag.
\newblock {\em Advanced mathematical methods for scientists and engineers}.
\newblock McGraw-Hill, New York, 1978.

\bibitem{Boyarsky(1997)}
Abraham Boyarsky and Pawel G\'ora.
\newblock {\em Laws of chaos : invariant measures and dynamical systems in one
  dimension}.
\newblock Birkh\"auser, Boston, Mass., 1997.

\bibitem{Calvo(1994)}
MP~Calvo, A.~Murua, and JM~Sanz-Serna.
\newblock Modified equations for {ODE}s.
\newblock In {\em Chaotic numerics: an International Workshop on the
  Approximation and Computation of Complicated Dynamical Behavior, Deakin
  University, Geelong, Australia, July 12-16, 1993}, volume 172 of {\em
  Contemporary Mathematics}, page~63. American Mathematical Society, 1994.

\bibitem{channell1990symplectic}
PJ~Channell and C.~Scovel.
\newblock Symplectic integration of {Hamiltonian} systems.
\newblock {\em Nonlinearity}, 3:231, 1990.

\bibitem{Corless(1992)continued}
Robert~M. Corless.
\newblock Continued fractions and chaos.
\newblock {\em The American Mathematical Monthly}, 99(3):203--215, 1992.

\bibitem{Corless(1994)}
Robert~M. Corless.
\newblock Error backward.
\newblock In P.~Kloeden and K.~Palmer, editors, {\em Proceedings of Chaotic
  Numerics, Geelong, 1993}, volume 172 of {\em AMS Contemporary Mathematics},
  pages 31--62, 1994.

\bibitem{Corless(1994)c}
Robert~M. Corless.
\newblock What good are numerical simulations of chaotic dynamical systems?
\newblock {\em Computers \& Mathematics with Applications}, 28(10-12):107--121,
  1994.

\bibitem{Corless(1992)}
Robert~M. Corless and George~F. Corliss.
\newblock Rationale for guaranteed {ODE} defect control.
\newblock In L.~Atanassova and J.~Herzberger, editors, {\em Computer Arithmetic
  and Enclosure Methods}, pages 3--12. North-Holland, 1992.

\bibitem{corless1991numerical}
Robert~M Corless, C~Essex, and MAH Nerenberg.
\newblock Numerical methods can suppress chaos.
\newblock {\em Physics Letters A}, 157(1):27--36, 1991.

\bibitem{corless2013graduate}
Robert~M. Corless and Nicolas Fillion.
\newblock {\em A graduate introduction to numerical methods}.
\newblock Springer-Verlag, 2013.

\bibitem{CorlessFillion2025SIAM}
Robert~M. Corless and Nicolas Fillion.
\newblock {\em Perturbation methods using backward error}.
\newblock SIAM, Philadelphia, 2025.

\bibitem{Corless2016}
Robert~M. Corless and Julia~E. Jankowski.
\newblock Variations on a theme of {E}uler.
\newblock {\em SIAM Review}, 58(4):775--792, January 2016.

\bibitem{Corless2018}
Robert~M. Corless, C.~Yal\c{c}in Kaya, and Robert H.~C. Moir.
\newblock Optimal residuals and the {D}ahlquist test problem.
\newblock {\em Numerical Algorithms}, 81(4):1253--1274, November 2018.

\bibitem{corless2019optimal}
Robert~M. Corless, C~Yal{\c{c}}{\i}n Kaya, and Robert~HC Moir.
\newblock Optimal residuals and the {Dahlquist} test problem.
\newblock {\em Numerical Algorithms}, 81(4):1253--1274, 2019.

\bibitem{crouch1993numerical}
P.E. Crouch and R.~Grossman.
\newblock Numerical integration of ordinary differential equations on
  manifolds.
\newblock {\em Journal of Nonlinear Science}, 3(1):1--33, 1993.

\bibitem{Enright(2007)}
Wayne~H. Enright and Wayne~B. Hayes.
\newblock Robust and reliable defect control for {Runge-Kutta} methods.
\newblock {\em ACM Trans. Math. Softw.}, 33, March 2007.

\bibitem{Griffiths(1986)}
D.F. Griffiths and J.M. Sanz-Serna.
\newblock On the scope of the method of modified equations.
\newblock {\em SIAM Journal on Scientific and Statistical Computing}, 7:994,
  1986.

\bibitem{hairer1997life}
E.~Hairer and C.~Lubich.
\newblock The life-span of backward error analysis for numerical integrators.
\newblock {\em Numerische Mathematik}, 76(4):441--462, 1997.

\bibitem{Hairer(2006)}
E.~Hairer, C.~Lubich, and G.~Wanner.
\newblock {\em Geometric numerical integration: structure-preserving algorithms
  for ordinary differential equations}.
\newblock Springer Verlag, 2006.

\bibitem{Hairer(1993)}
Ernst Hairer, Syvert~P. N{\o}rsett, and Gerhard Wanner.
\newblock {\em Solving ordinary differential equations: Nonstiff problems}.
\newblock Springer, 1993.

\bibitem{hayes2007fast}
Wayne~B Hayes and Kenneth~R Jackson.
\newblock A fast shadowing algorithm for high-dimensional {ODE} systems.
\newblock {\em SIAM Journal on Scientific Computing}, 29(4):1738--1758, 2007.

\bibitem{kunkel2006differential}
Peter Kunkel and Volker Mehrmann.
\newblock {\em Differential-algebraic equations: analysis and numerical
  solution}, volume~2.
\newblock European Mathematical Society, 2006.

\bibitem{kunkel2023discretization}
Peter Kunkel and Volker Mehrmann.
\newblock Discretization of inherent {ODE}s and the geometric integration of
  {DAE}s with symmetries.
\newblock {\em BIT Numerical Mathematics}, 63(2):29, 2023.

\bibitem{liao2006chaos}
Xiaoxin Liao and Pei Yu.
\newblock Chaos control for the family of {R}{\"o}ssler systems using feedback
  controllers.
\newblock {\em Chaos, Solitons \& Fractals}, 29(1):91--107, 2006.

\bibitem{lorenz1963deterministic}
Edward~N Lorenz.
\newblock Deterministic nonperiodic flow.
\newblock {\em Journal of atmospheric sciences}, 20(2):130--141, 1963.

\bibitem{Pilyugin(1999)}
Sergei~Yu Pilyugin.
\newblock {\em Shadowing in dynamical systems}, volume 1706.
\newblock Springer, Berlin ; New York, 1999.

\bibitem{Preto(1999)class}
M.~Preto and S.~Tremaine.
\newblock A class of symplectic integrators with adaptive time step for
  separable {Hamiltonian} systems.
\newblock {\em The Astronomical Journal}, 118:2532, 1999.

\bibitem{rackauckas2017differentialequations}
Christopher Rackauckas and Qing Nie.
\newblock Differentialequations. jl--a performant and feature-rich ecosystem
  for solving differential equations in {Julia}.
\newblock {\em Journal of open research software}, 5(1):15--15, 2017.

\bibitem{SanzSerna(1994)}
J.M. Sanz-Serna and M.P. Calvo.
\newblock {\em Numerical {H}amiltonian problems}.
\newblock Chapman \& Hall/CRC, 1994.

\bibitem{Sparrow:1982}
Colin Sparrow.
\newblock {\em The Lorenz equations :bifurcations, chaos, and strange
  attractors}, volume~41.
\newblock Springer-Verlag, New York, 1982.

\bibitem{stuart1998dynamical}
Andrew Stuart and Anthony~R Humphries.
\newblock {\em Dynamical systems and numerical analysis}, volume~2.
\newblock Cambridge University Press, 1998.

\bibitem{Tupper2009}
Paul Tupper.
\newblock The relation between approximation in distribution and shadowing in
  molecular dynamics.
\newblock {\em SIAM Journal on Applied Dynamical Systems}, 8(2):734--755, 2009.

\end{thebibliography}

\appendix
\section{Julia script for half-precision simulation of the Gauss map}
\begin{verbatim}
#
# Copyright (c) 2024 Robert M. Corless
#
# Permission is hereby granted, free of charge, to any person obtaining a copy
# of this software and associated documentation files (the "Software"), to
# deal in the Software without restriction, including without limitation the
# rights to use, copy, modify, merge, publish, distribute, sublicense, and/or
# sell copies of the Software, and to permit persons to whom the Software is
# furnished to do so, subject to the following conditions:
#
# The above copyright notice and this permission notice shall be included in
# all copies or substantial portions of the Software.
#
# THE SOFTWARE IS PROVIDED "AS IS", WITHOUT WARRANTY OF ANY KIND, EXPRESS
# OR IMPLIED, INCLUDING BUT NOT LIMITED TO THE WARRANTIES OF MERCHANTABILITY,
# FITNESS FOR A PARTICULAR PURPOSE AND NONINFRINGEMENT. IN NO EVENT SHALL
# THE AUTHORS OR COPYRIGHT HOLDERS BE LIABLE FOR ANY CLAIM, DAMAGES OR OTHER
# LIABILITY, WHETHER IN AN ACTION OF CONTRACT, TORT OR OTHERWISE, ARISING
# FROM, OUT OF OR IN CONNECTION WITH THE SOFTWARE OR THE USE OR OTHER
# DEALINGS IN THE SOFTWARE.

function frac(x)
  return x - floor(x)
end

# Gauss map from [0,1] to [0,1]
function G(x)
  if x==0
    return 0
  else
    return frac(1/x)
  end
end

# The number of float16s in [0,1] inclusive
# is 15*1024 + 1 = 15361
N = 15361
X = Array{Float16}(undef, N)
Y = Array{Float16}(undef, N)

# Start the iteration exactly at 1 (G(1) will be zero)
X[1] = Float16(1)
Y[1] = G(X[1])
# Step through all the float16s and see what G does to each
for k=2:N
  X[k] = prevfloat(X[k-1])
  Y[k] = G(X[k])
end

# Recording is a bit of a problem.  We do not wish to
# use the decimal versions, because of binary-to-decimal-to-binary
# conversion questions.  So we record instead the indices into X.

using CSV, Tables

# X is sorted in decreasing order
# E will be the table of edges in the graph of G: Float16 --> Float16.

E = Array{Int16}(undef,N)
for j=1:N
  # Look for Y[j] in X
  y = Y[j]
  k = 1
  # NaN16 fails this test so k=1 for NaN16s
  while X[k]>y && k<N
    k = k+1
  end
  E[j] = k
end

# X[1] = Float16(1)
# X[k] = prevfloat(X[k-1]) for 1 < k <= 15361
# Last X to produce a nonzero G is X[10224] = 0.000993
# Last X to produce a non-NaN G is X[15104] = 1.53e-5
# First X with G(X)=Nan16 is X[15105] = 1.526e-5
# X[15361] = Float16(0)

CSV.write( "edges.csv", Tables.table( E ) )
\end{verbatim}
\end{document}